\documentclass[journal]{IEEEtran}

\usepackage{graphicx}
\usepackage{caption}

\usepackage{algorithm}
\usepackage[noend]{algpseudocode}

\usepackage{multicol}
\usepackage{amsmath,amssymb}
\usepackage{amsfonts}
\usepackage{textcomp}
\usepackage{psfrag}
\usepackage{multimedia}
\usepackage{enumitem}

\newtheorem{theorem}{Theorem}[section]

\newtheorem{Lemma}[theorem]{Lemma}
\newtheorem{Definition}[theorem]{Definition}
\newtheorem{Remark}[theorem]{Remark}
\newtheorem{Example}[theorem]{Example}
\newtheorem{Assumption}[theorem]{Assumption}

\usepackage[usenames,dvipsnames]{color}
\definecolor{wheat}{rgb}{0.96,0.87,0.70}
\definecolor{mario}{rgb}{0.8,0.8,1}
\definecolor{seb}{rgb}{0.8,1,0.8}
\definecolor{myGreen}{rgb}{0.,0.8,0.0}

\definecolor{darkgreen}{rgb}{0,0.6,0}

\def\R{\mathbb{R}}

\newcommand{\ontop}[2]{\genfrac{}{}{0pt}{}{#1}{#2}}

\newcounter{lastnote}

\usepackage{tikz}
\usetikzlibrary{calc,arrows,positioning}

\begin{document} 

\title{Stabilization of strictly pre-dissipative nonlinear receding horizon control by terminal costs\thanks{The research for this paper was supported by DFG Grant No.\ 499435839.}}

\author
{Lars Gr\"une, Mario Zanon
	\thanks{Lars  Gr\"une is with the University of Bayreuth, Germany. E-Mail: lars.gruene@uni-bayreuth.de}
	\thanks{Mario Zanon is with the IMT School for Advanced Studies Lucca, Italy. }
	\thanks{}
}

\IEEEtitleabstractindextext{
	\begin{abstract}
	It is known that receding horizon control with a strictly pre-dissipative optimal control problem yields a practically asymptotically stable closed loop when suitable state constraints are imposed. In this note we show that alternatively suitably bounded terminal costs can be used for stabilizing the closed loop.
	\end{abstract}
	
}

\maketitle

\IEEEdisplaynontitleabstractindextext

\IEEEpeerreviewmaketitle

\section{Introduction}

Receding horizon control (often used synonymously with model predictive control) is a control technique in which a finite horizon optimal control problem is solved in each time step and the first element of the resulting optimal control sequence is used in the next time step \cite{Rawlings2017,Gruene2017}. Under suitable stabilizability and regularity conditions, this scheme yields a practically asymptotically stable closed loop if the system is strictly dissipative with supply function defined via the stage cost of the finite horizon optimal control problem \cite[Chapter 8]{Gruene2017}. In this case, we call the optimal control problem {\em strictly dissipative}. Here, the size of the ``practical'' neighborhood of the equilibrium to which the closed-loop solution converges is determined by the length of the finite optimization horizon. True (as opposed to practical) asymptotic stability can be achieved by using suitable terminal constraints and costs, see \cite{DiAR10,Amrit2011a} or Theorem 8.13 in \cite{Gruene2017}. In these approaches the terminal cost is typically a local control Lyapunov function for the system and the terminal constraints are needed because the design of a global control Lyapunov function is usually a very difficult task. As a simpler alternative, it was shown in \cite{Zanon2018a} that linear terminal costs can also be used to obtain true asymptotic stability.

The strict dissipativity property that is at the heart of all these results requires the existence of a so-called storage function $\lambda$ mapping the state space into the reals. It is a strengthened version of the system theoretic dissipativity property introduced by Willems in his seminal papers \cite{Willems1972a,Willems1972b} and also featured in his slightly earlier paper \cite{Willems1971} on linear quadratic optimal control and the algebraic Riccati equation. Readers familiar with Lyapunov's stability theory can see the storage function $\lambda$ as a generalization of a Lyapunov function. However, unlike Lyapunov functions, $\lambda$ need not attain nonnegative values. However, it must be bounded from below, and this property is crucial for deriving the (practical) stability properties for receding horizon control cited above.

For generalized linear quadratic problems, i.e., problems with linear dynamics and a cost function containing quadratic and linear terms, with state space $\R^{n_x}$, a standard construction for a storage function results in a function of the form $\lambda(x)= x^TPx + \nu^Tx$, for $P\in\R^{n_x\times n_x}$ and $\nu\in\R^{n_x}$, see \cite[Proposition 4.5]{DGSW14}. Clearly, such a function $\lambda$ is in general {\em not} bounded from below and Example 2.3 in \cite{DGSW14}, which we also present as Example \ref{ex:simple_lqr1}, below, shows that storage functions unbounded from below may occur even for very simple scalar problems. While the potential unboundedness of $\lambda$ has been handled somewhat informally in \cite{DGSW14}, later in \cite{Gruene2018} the variant of strict dissipativity with storage function not bounded from below has been termed {\em strict pre-dissipativity}. For strictly pre-dissipative problems, one way to obtain strict dissipativity and thus (practical) asymptotic stability is to suitably restrict the state space by means of state constraints, e.g., to a compact set, on which $\lambda$ is bounded from below. 

Since such a restriction of the state space may not always be desirable, in this note we will look at an alternative way to regain (practical) asymptotic stability. More precisely, we want to answer the following question: Given a receding horizon control scheme with strictly pre-dissipative optimal control problem, can we add a simple terminal cost that guarantees (practical) asymptotic stability? Here ``simple'' means that we don't want to design control Lyapunov function terminal costs but terms that are easier to compute. We will see that for obtaining practical asymptotic stability is is sufficient that the terminal cost is larger than the negative storage function, while for obtaining asymptotic stability we need that the sum of the storage function and the terminal cost is positive semidefinite at the optimal equilibrium. It is worth to be noting that this implies the necessary condition from \cite{Zanon2018a}, cf.\ the discussion after Assumption \ref{as:semidef}, below. We emphasize that no terminal constraints are needed in our approach. This paper focuses on general nonlinear optimal control problems. In the companion paper \cite{ZanG24} we discuss the linear-quadratic case, for which further results can be obtained.

The remainder of this note is organized as follows. In Section \ref{sec:problem} we define the problem and the concepts we use. Section \ref{sec:results} contains the main results and proofs, in which we heavily rely on \cite{GMTT05} and \cite[Chapter 8]{Gruene2017}. Section \ref{sec:example} contains an illustrative example and Section \ref{sec:conclusion} concludes the note.

\section{Problem Statement}\label{sec:problem}

\subsection{Receding Horizon Control}

We consider discrete-time systems of the form
\begin{align}
	\label{eq:system}
	x_{k+1} = f(x_k,u_k),
\end{align}
where $x\in\mathbb{R}^{n_x}$ and $u\in\mathbb{R}^{n_u}$ denote the states and controls respectively. 

Receding horizon or Model predictive control consists in minimizing a given stage cost $\ell:\R^{n_x}\times\R^{n_u}\to\R$ over a fixed finite prediction horizon $N$, possibly subject to constraints and with the addition of a terminal cost. The receding horizon optimal control problem (RH-OCP) reads
\begin{subequations}
	\label{eq:empc}
	\begin{align}
		\min_{x_0,u_0,\ldots,x_N} \ \ & \sum_{k=0}^{N-1} \ell(x_k,u_k) + V^\mathrm{f}(x_N) \label{eq:empc_cost}\\
		\mathrm{s.t.} \ \ & x_0 = \hat x_j, \label{eq:empc_ic}\\
		&x_{k+1} = f(x_k,u_k),  && k \in\mathbb{I}_0^{N-1}, \label{eq:empc_dyn}\\
		&h(x_k,u_k) \leq 0,  && k \in\mathbb{I}_0^{N-1}, \label{eq:empc_pc}
	\end{align}
\end{subequations}
where $h:\R^{n_x}\times\R^{n_u}\to\R^l$ defines the state and input constraints and inequality \eqref{eq:empc_pc} is to be understood componentwise. We assume that all involved functions $f,\ell,V^\mathrm{f},h$ are continuous. It is possible to introduce additional constraints on the terminal state $x_N$, but we will not consider that option in this note. The solution to this problem, which we assume to exist for all $N\in\mathbb{N}$ and all $x_0\in\R^n$ satisfying the constraints \eqref{eq:empc_pc} for some $u_0\in\R^{n_u}$, is denoted by $x_k^\star$, $u_k^\star$.

We refer to this optimal solution as the MPC prediction. 
In order to distinguish the time in the MPC closed loop from the time $k$ of the MPC prediction, we denote the former by $j$.
Starting in the initial value $\hat x_0$ at time instant $j=0$, for every time instant $j\ge 0$ the state $\hat x_j$ is measured, Problem~\eqref{eq:empc} is solved with $x_0=\hat x_j$, and the first optimal input $u_0^\star$ is applied to the system to obtain
\begin{align}
	\label{eq:clsystem}
	 \hat x_{j+1} = f(\hat x_j, u_0^\star). 
\end{align}%
This procedure is repeated iteratively for all $j\ge 0$.

Associated to the RH-OCP~\eqref{eq:empc}, we define the steady-state optimization problem (SOP) 
\begin{align}
	\label{eq:sop_dt}
	\min_{\bar x, \bar u} \ \ & \ell(\bar x, \bar u) &
	\mathrm{s.t.} \ \  
	& \bar x = f(\bar x,\bar u), &
	& h(\bar x, \bar u)\leq 0.
\end{align}
The optimal solution of SOP~\eqref{eq:sop_dt} is denoted as $\bar z^\star = (\bar x^\star,\bar u^\star)$. 

In this note, we are interested in obtaining stability properties of the closed loop system~\eqref{eq:clsystem} at the optimal steady state $\bar z^\star$. While stability results are abundant for the case of suitably formulated terminal constraints and Lyapunov function terminal costs~\cite{Diehl2011, Angeli2012a,Faulwasser2018a,Muller2015,Muller2013a}, we focus here on the case of no terminal constraints. This case has been analyzed, e.g., in~\cite{Faulwasser2018,Zanon2018a,Grune2013a,Muller2016}, where we can further distinguish between formulations without terminal cost and formulations with simple terminal costs that need not be Lyapunov functions, which are usually difficult to design. This last approach has in particular been taken in~\cite{Faulwasser2018,Zanon2018a} by using a linear terminal cost and the present note can be seen as a continuation of this research. As in these references, our analysis is based on dissipativity concepts and the turnpike property.

\subsection{Strict dissipativity}

The stability theory of receding horizon control and economic MPC is often based on strict dissipativity~\cite{Diehl2011, Angeli2012a,Faulwasser2018a,Muller2015,Muller2013a,Grune2013a}. Next we define the weaker notion of {\em strict pre-dissipativity}, which we will use throughout this note.
\begin{Definition}
	We say that the RH-OCP \eqref{eq:empc} is {\em strictly pre-dissipative} if there exists a continuous {\em storage function} $\lambda:\R^{n_x}\to\R$ such that for all $(x,u)\in\mathcal{Z}$ the {\em rotated cost}
	\begin{align}
		\label{eq:rotated cost}
		L(x,u):=\ell(x,u) - \ell(\bar x^\star,\bar u^\star) + \lambda(x) - \lambda(f(x,u))
	\end{align}
	satisfies
	\begin{align}
		\label{eq:strict_dissipativity}
		L(x,u) \geq \alpha (\|x\|),
	\end{align}
	where $\alpha$ is a class $\mathcal{K}$ function.
\label{def:diss}\end{Definition}
In contrast to strict dissipativity, strict pre-dissipativity, introduced under this name in \cite{Gruene2018}, does not require the storage function $\lambda$ to be bounded from below. This implies that one cannot use arguments as, e.g., in~\cite{Grune2013a,Gruene2014} in order to conclude (practical) stability properties of the closed loop \eqref{eq:clsystem}, and in fact stability may fail to hold, as we will show by means of the following example.

\begin{Example}\label{ex:simple_lqr1}
Consider the optimal control problem with dynamics and stage cost
\[ x_{k+1} = 2x_k+u_k, \qquad \ell(x,u) = u^2.\]
One easily sees that for any initial condition $x_0$ and any horizon $N$ the optimal control sequence is $u_k^\star \equiv 0$, as this is the only control that produces $0$ cost, while all other control sequences produce positive costs. This implies that system \eqref{eq:clsystem} becomes 
\[ \hat x_{j+1} = 2 \hat x_j, \]
for which the origin is obviously exponentially unstable. Yet, one checks that this problem is strictly pre-dissipative at the optimal equilibrium $\bar z^\star = (0,0)$ with storage function $\lambda(x) = - c x^2$ for each $c\in(0,1]$. This shows that strict pre-dissipativity does not imply asymptotic stability of the optimal equilibrium.
\end{Example}

As already mentioned in the introduction and as also seen in this example, storage functions that are not bounded from below appear naturally already for linear quadratic problems. In order to achieve closed-loop stability, often a compact state constraint set is imposed, as compactness implies boundedness of the storage function provided it is continuous (which is often the case). For Example \ref{ex:simple_lqr1}, it was shown in \cite[Example 2.3]{DGSW14} that this indeed renders the origin practically asymptotically stable for the closed loop. Yet, imposing compact state constraints just for the sake of achieving stability may not always be desirable.
As we will prove in this note, stability can be alternatively achieved by a suitably defined terminal cost, cf.\ the end of Section \ref{sec:lq}, below.

\section{Main results}\label{sec:results}

We start with the following lemma, which shows the relation between the optimal control problems for the original and the rotated stage cost, respectively.

\begin{Lemma}
	\label{lem:rotated_equivalence_nonlinear}
Consider the RH-OCP \eqref{eq:empc} for an arbitrary finite horizon $N$. Assume strict pre-dissipativity with storage function $\lambda$. Then the problem with stage cost $\ell$ and terminal cost $V^{\mathrm{f}}$ has the same optimal trajectories $x_k^\star$ and control sequences $u_k^\star$ as the RH-OCP problem with rotated stage cost $L$ from \eqref{eq:rotated cost} and adapted terminal cost $V^{\mathrm{f}}+\lambda$. Moreover, the corresponding optimal value functions satisfy $V_{\lambda,N} = V_N + \lambda- N\ell(\bar x^\star,\bar u^\star)$.
\end{Lemma}
\begin{IEEEproof}
	The proof follows from the observation that
	\begin{align*}
	  &\sum_{k=0}^{N-1} L(x_k,u_k) + V^{\mathrm{f}}(x_N)+\lambda(x_N) \\
	  & = \sum_{k=0}^{N-1} \big[\ell(x_k,u_k) - \ell(\bar x^\star,\bar u^\star) + \lambda(x_k) - \lambda(x_{k+1})\big]\\
	  & \qquad + V^{\mathrm{f}}(x_N)+\lambda(x_N)\\
	  & = \lambda(x_0) - N\ell(\bar x^\star,\bar u^\star) + \sum_{k=0}^{N-1} \ell(x_k,u_k) + V^{\mathrm{f}}(x_N)
	\end{align*}
	implies that the optimization objectives of the two problems differ only by the constant $\lambda(x_0)- N\ell(\bar x^\star,\bar u^\star)$. Note that all values in the above sums are finite, as $\lambda(x)\in\R$ for all $x\in\R^{n_x}$. From this, both the statements about the optimal solutions, i.e., the minimizers as well as the statement about the optimal value functions follow.
\end{IEEEproof}

Lemma \ref{lem:rotated_equivalence_nonlinear} now allows us to use existing results on stability of MPC for either stabilizing or strictily dissipative optimal control problems and carry them over to the strictly pre-dissipative case. For doing this, we distinguish between two cases.

\subsection{Results for positive semidefinite $V^{\mathrm{f}}+\lambda$}

In this subsection we make the assumption that $V^{\mathrm{f}}+\lambda$ is positive semidefinite in the following sense. This will allow us to use stability results from \cite{GMTT05}. 

\begin{Definition}
	A function $\Phi:\R^{n_x}\to\R$ is called positive semidefinite at a point $x^\star\in\R^{n_x}$, if $\Phi(x^*)=0$ and $\Phi(x)\ge 0$ for all $x\in\R^{n_x}$. 
\end{Definition}

The main structural assumption we make in this subsection is the following, where we use the state and control constraint sets  
\[ \mathcal{X} := \{ x\in\R^{n_x}\,|\, \mbox{there is } u \in \R^{n_u} \mbox{ with } h(x,u)\le 0\}, \]
and
\[ \mathcal{U} := \{ u\in\R^{n_u}\,|\, \mbox{there is } x \in \R^{n_x} \mbox{ with } h(x,u)\le 0\}. \]

\begin{Assumption}\label{as:semidef}
	The optimal control problem \eqref{eq:empc} is strictly pre-dissipative at an equilibrium $x^\star\in\mathrm{int}\mathcal{X}$ with storage function $\lambda$ and class $\mathcal{K}$ function $\alpha$ and the function $x\mapsto V^{\mathrm{f}}(x)+\lambda(x)$ is positive semidefinite at $x^\star$.
\end{Assumption}

We note that this assumption implies that $V^{\mathrm{f}}+\lambda$ has a global minimum at $x^\star$, which satisfies the necessary optimality conditions of a local minimum since is is in the interior of $\mathcal{X}$. If both functions are differentiable, this implies that $\nabla V^{\mathrm{f}}(x^\star) = - \nabla \lambda(x^\star)$ must hold. Hence, the linear term in the terminal cost $V^{\mathrm{f}}$ provides a gradient correction in the sense of \cite{Faulwasser2018,Zanon2018a}, which is known to be necessary for obtaining asymptotic stability in the absence of terminal constraints.

In addition, for invoking the results from \cite{GMTT05} we need the following technical assumption on the rotated stage cost.

\begin{Assumption} \label{as:technical} 
Either $\mathcal{U}$ is compact or for each compact set $C\subset \mathcal{X}$, each $N\ge 1$, and each $\mu>0$ there is $\eta>0$ such that 
\[ \sum_{k=0}^{N-1} L(x_k,u_k) + V^{\mathrm{f}}(x_N)+\lambda(x_N) \le \eta\]
implies $\|u_k\|\le \mu$ for all $k=0,\ldots,N-1$.
\end{Assumption}
We note that the second alternative of the assumption follows, e.g., if $\ell(x,u) = \ell_1(x) + \ell_2(u)$ and $\ell_2(u) \ge \gamma(\|u\|)$ for some $\gamma\in\mathcal{K}$ and $V^{\mathrm{f}}$ is bounded from below.

\begin{theorem} Consider the MPC closed loop \eqref{eq:clsystem} with optimal control problem \eqref{eq:empc} satisfying Assumptions \ref{as:semidef} and \ref{as:technical}. Assume there is $\rho\in\mathcal{K}$ such that the optimal value function of \eqref{eq:empc} satisfies 
\[ V_N(x) +\lambda(x) - N\ell(\bar x^\star, \bar u^\star) \le \inf_{u\in \mathcal{U}} \rho(L(x,u)) \] 
for all $x\in\mathcal{X}$ and all $N\ge 1$. Then there are $\beta\in\mathcal{KL}$, and $\Delta(N)>\delta(N)>0$ with $\Delta(N)\to \infty$ and $\delta(N)\to 0$ as $N\to\infty$, such that for all sufficiently large $N$ the solutions $\hat x_j$ of the closed loop \eqref{eq:clsystem} with $\|\hat x_0\| \le \Delta(N)$ satisfy
\begin{equation} \|\hat x_j\| \le \max\{\beta(\|\hat x_0\|,j),\, \delta(N)\}.\label{eq:spas}\end{equation}
In words, $x^\star$ is a  semiglobally practically asymptotically stable equilibrium of the closed loop \eqref{eq:clsystem}.

If the inequality for $V_N$ holds with a linear function $\rho\in\mathcal{K}$, then $x^\star$ is an asymptotically stable equilibrium of the closed loop \eqref{eq:clsystem}, i.e., inequality \eqref{eq:spas} holds with $\delta(N)=0$ for all initial conditions $\hat x_0$.
\label{thm:semidef}
\end{theorem}
\begin{IEEEproof} For the problem with rotated stage cost and adapted terminal cost, semiglobal practical asymptotic stability follows from \cite[Corollary 1]{GMTT05} and ``true'' asymptotic stability in case of $\rho$ being linear follows from \cite[Corollary 3]{GMTT05}. Since by Lemma \ref{lem:rotated_equivalence_nonlinear} the optimal controls for the original problem coincide with the optimal controls for the rotated problem, the resulting closed loop systems are identical, implying the result for the original problem \eqref{eq:empc}.

It remains to verify the assumptions of Corollary 1 and 3 in \cite{GMTT05} for the rotated problem. The standing assumption SA1 in \cite{GMTT05} is assumed right after \eqref{eq:empc}. SA2 follows from Assumption \ref{as:technical}. SA3 follows, using Remark 1 in \cite{GMTT05}, with $\sigma=\alpha$ and $W\equiv 0$ from the positive definiteness of the rotated stage cost $L$, which in turn follows from strict pre-dissipativity. Finally, SA4 is contained in the assumptions of Theorem \ref{thm:semidef}, with $\alpha_W = \rho$, $\gamma_W={\rm Id}$ and $\bar\gamma_W$ arbitrary. This shows that all the assumptions for \cite[Corollary 1]{GMTT05} hold. The additional assumption needed for Corollary 3 in this reference requires $\alpha_W, \gamma_W$, and $\bar\gamma_W$ to be linear functions, which is satisfied if $\rho$ is linear. Hence, both corollaries can be used and show the desired stability properties.
\end{IEEEproof}
\begin{Remark} (i) We note that \cite[Corollary 3]{GMTT05} also provides an estimate for the horizon length $N$ for which asymptotic stability holds in case of a linear $\rho$, but this bound is very conservative. Tighter bounds were provided, e.g., in \cite{TuMT06} or \cite{GPSW10}, see also \cite[Chapter 6]{Gruene2017}.

(ii) The condition on $V_N$ in Theorem \ref{thm:semidef} is effectively a condition on the optimal value function $V_{\lambda,N}$ of the problem with rotated cost and adapted terminal cost, cf.\ Lemma \ref{lem:rotated_equivalence_nonlinear}. In essence, it requires that the system can be controlled asymptotically to $\bar x^\star$ with sufficiently low cost. In particular, if $L$ and the terminal cost terms are polynomial and the system can be controlled to $\bar x^\star$ exponentially fast, then $\rho$ can be chosen as a linear function. We refer to \cite[Chapter 6]{Gruene2017} for an extensive discussion and examples.
\label{rem:semidef}\end{Remark}

\subsection{Results for $V^{\mathrm{f}}+\lambda$ bounded from below}

The design of a terminal cost $V^\mathrm{f}$ for which $V^{\mathrm{f}}+\lambda$ is positive semidefinite requires rather accurate knowledge of the storage function $\lambda$. This may either not be available or difficult to obtain. To this end, we now present a semiglobal practical asymptotic stability result under the following significantly weaker assumption.

\begin{Assumption}\label{as:bounded}
The optimal control problem \eqref{eq:empc} is strictly pre-dissipative at an equilibrium $x^\star\in\mathrm{int}\mathcal{X}$ with storage function $\lambda$ and class $\mathcal{K}$ function $\alpha$ and the function $x\mapsto V^{\mathrm{f}}(x)+\lambda(x)$ is bounded from below on $\mathcal{X}$.
\end{Assumption}
In contrast to requiring positive semidefiniteness, this lower boundedness assumption only requires $V^{\mathrm{f}}$ to be sufficiently large for large $x$, as on compact sets $V^{\mathrm{f}}+\lambda$ is bounded from below by continuity, regardless of how $V^{\mathrm{f}}$ is chosen.

Again, besides this structural assumption we need a couple of technical assumptions. Here we use the assumptions from \cite[Chapter 8]{Gruene2017}, which are somewhat more streamlined than the assumptions in the original paper \cite{Gruene2014}, where the result we use appeared for the first time.

\begin{Assumption}\label{as:technical2}
There are class $\mathcal{K}$ functions $\gamma_{V}$ and $\gamma_{V_{\lambda}}$, such that the inequalities
\[ |V_N(x)-V_N(\bar x^\star)| \le \gamma_V(\|x-\bar x^\star\|) \]
and
\[ |V_{\lambda,N}(x)-V_{\lambda,N}(\bar x^\star)| \le \gamma_{V_\lambda}(\|x-\bar x^\star\|) \]
hold for all $N\ge 1$ and all $x\in\mathcal{X}$.
\end{Assumption}

\begin{theorem} Consider the MPC closed loop \eqref{eq:clsystem} with optimal control problem \eqref{eq:empc} satisfying Assumptions \ref{as:bounded} and \ref{as:technical2}. Assume there is $\rho\in\mathcal{K}$ such that the optimal value function of \eqref{eq:empc} satisfies 
\[ V_N(x) +\lambda(x) - N\ell(\bar x^\star, \bar u^\star) \le \inf_{u\in \mathcal{U}} \rho(L(x,u)) \] 
for all $x\in\mathcal{X}$ and all $N\ge 1$. Then there are $\beta\in\mathcal{KL}$, and $\Delta(N)>\delta(N)>0$ with $\Delta(N)\to \infty$ and $\delta(N)\to 0$ as $N\to\infty$, such that for all sufficiently large $N$ the solutions $\hat x_j$ of the closed loop \eqref{eq:clsystem} with $\|\hat x_0\| \le \Delta(N)$ satisfy
\begin{equation} \|\hat x_j\| \le \max\{\beta(\|\hat x_0\|,j),\, \delta(N)\}.\label{eq:spas2}\end{equation}
In words, $x^\star$ is a semiglobally practically asymptotically stable equilibrium of the closed loop \eqref{eq:clsystem}.
\label{thm:bounded}\end{theorem}
\begin{IEEEproof}
The proof proceeds almost identical to the proof of \cite[Theorem 8.33]{Gruene2017}, with the following two changes:
\begin{itemize}
\item In the proof of \cite[Theorem 8.33]{Gruene2017} practical asymptotic stability is obtained by showing that the optimal value function for the problem with stage cost $L$ and without terminal cost is a practical Lyapunov function for the closed loop generated by the problem with stage cost $L$ and terminal cost $\lambda$. As no particular storage function properties of $\lambda$ are exploited in this proof, we can apply the same reasoning to the closed loop generated by the problem with stage cost $L$ and terminal cost $V^\mathrm{f}+\lambda$ instead of $\lambda$, which by Lemma \ref{lem:rotated_equivalence_nonlinear} yields the same closed loop system as \eqref{eq:empc}; see also \cite{Faulwasser2018} for a similar reasoning.
\item In \cite[Theorem 8.33]{Gruene2017}, boundedness of $\lambda$ (which, as just mentioned, needs to be replaced here by $V^{\mathrm{f}}+\lambda$) is required. However, an inspection of the proof shows that the upper bound on this function is only used when the function is evaluated in the initial condition. Hence, for any fixed $\Delta>0$ we can apply the proof for all solutions with initial conditions $\|\hat x_0\|\le \Delta$. Thus we obtain a lower bound on $N$, depending on $\Delta$, for which practical asymptotic stability holds for these solutions. As we can find such bounds on $N$ for arbitrary $\Delta>0$, this implies the claimed semiglobal practical stability.
\end{itemize}
\end{IEEEproof}

\subsection{Relation to the required supply}

One of the standard constructions of a storage function for a dissipative system with a given supply rate is the so-called required supply. Here we follow the definition in \cite{LopS06}. In order to adapt this definition to our strict setting, we assume that the optimal control problem \eqref{eq:empc} satisfies Definition \ref{def:diss} and define the supply rate via $s(x,u) := \ell(x,u) - \ell(\bar x^\star,\bar u^\star) - \alpha(\|x-\bar x^\star\|)$, with $\alpha$ from Definition \ref{def:diss}. For the definition of the required supply we need to make the assumption that every $x\in\R^{n_x}$ can be reached from $\bar x^\star$. Then, the function defined by 
\[ \lambda_{rs}(x) := \inf_{\ontop{u_0,\ldots,u_k,N_x}{x_0=\bar x^\star, x_{N_x} = x}} \sum_{k=0}^{N_x-1} s(x_k,u_k) \]
is finite and it is a storage function for the strict dissipativity property in Definition \ref{def:diss}, cf.\ \cite[Theorem 3.2]{LopS06}. 

As shown in \cite[Theorem 2(ii)]{Willems1972a}, the required supply is the largest possible storage function satisfying $\lambda_{rs}(\bar x^\star)=0$. This means that by choosing $\lambda = \lambda_{rs}$ in Theorems \ref{thm:semidef} and \ref{thm:bounded}, we obtain the least demanding condition on $V^{\mathrm{f}}$. 

Moreover, the choice $\lambda = \lambda_{rs}$ also gives an intuitive explanation for the requirements in Assumptions \ref{as:semidef} and \ref{as:bounded}. First note that Assumption \ref{as:semidef} implies $V^{\mathrm{f}} \ge - \lambda_{rs}$ while under Assumption \ref{as:bounded} we can assume $V^{\mathrm{f}} \ge - \lambda_{rs}$ without loss of generality, because adding a constant to $V^{\mathrm{f}}$ does not change the solutions to the optimal control problem. Moreover, even if one of the assumptions (or both) are violated, we can always without loss of generality assume that $V^{\mathrm{f}}(\bar x^\star) \ge - \lambda_{rs}(\bar x^\star)=0$, again because adding a constant to $V^{\mathrm{f}}$ does not change the optimal solutions. 

We provide the intuitive explanation of the two assumptions for such a choice of $V^{\mathrm{f}}$. To this end, consider an arbitrary solution $x_k$ with control $u_k$ starting in $x_0=\bar x^\star$ and reaching $x_N\ne \bar x^\star$ in $N\ge 1$ steps. Then, by nonnegativity of $\alpha(r)$, the definition of $\lambda_{rs}$ yields that 
\begin{align*} \lambda_{rs}(x) & \le \sum_{k=0}^{N-1} s(x_k,u_k)\\
& = \sum_{k=0}^{N-1} \Big( \ell(x_k,u_k) - \ell(\bar x^\star,\bar u^\star) - \alpha(\|x_k-\bar x^\star\|\Big) \\
& \le \sum_{k=0}^{N-1} \ell(x_k,u_k) - N\ell(\bar x^\star,\bar u^\star).
\end{align*}
Now, if we assume that $V^{\mathrm{f}} \ge - \lambda_{rs}$ is violated in $x_N$, i.e., if $V^{\mathrm{f}}(x_N) < -\lambda_{rs}(x_N)$, then the overall cost of the trajectory $x_k$ satisfies
\begin{align*} \sum_{k=0}^{N-1} \ell(x_k,u_k) + V^{\mathrm{f}}(x_N) & < \sum_{k=0}^{N-1} \ell(x_k,u_k) - \lambda_{rs}(x_N)\\
&\le N\ell(\bar x^\star,\bar u^\star)\\ 
&\le N\ell(\bar x^\star,\bar u^\star) + V^{\mathrm{f}}(\bar x^\star).\end{align*}
Observing that the last expression is precisely the cost to stay in $\bar x^\star$ for $N$ steps, this means that the finite horizon optimal solutions will not stay in $\bar x^\star$ because it is cheaper to move to $x_N$. While this does not exactly contradict stability of $\bar x^\star$ for the RH closed loop (the solutions could still stay in or near $\bar x^\star$ for a couple of steps steps before moving towards $x_N$), it is a clear indication that the optimal control problem is not well designed for obtaining asymptotic stability in or near $\bar x^\star$. 

From this observation, one may conjecture that $V^{\mathrm{f}} \ge - \lambda_{rs}$ is a necessary condition for obtaining (practical) asymptotic stability of the RH closed loop, but a formal proof appears technically involved and is beyond the scope of this note. However, in a companion paper we will present a tight lower bound for linear-quadratic problems.

\subsection{Consequences for the linear quadratic case}\label{sec:lq}

In the particular case of strictly dissipative generalized linear quadratic optimal control problems, it was shown in~\cite{Gruene2018} that the storage function can always be chosen to be of the form $\lambda(x) = x^\top\Lambda x + v^\top x$. Moreover, the technical Assumption \ref{as:technical} is always satisfied if $R>0$ and the technical Assumption \ref{as:technical2} as well as the bound on $V_N$ assumed in Theorem \ref{thm:semidef} are satisfied if $(A,B)$ are stabilizable. Hence, for applying Theorem \ref{thm:semidef} and Theorem \ref{thm:bounded}, it suffices to check the Assumptions \ref{as:semidef} and \ref{as:bounded}, respectively. 

If we restrict ourselves to terminal cost functions that contain only quadratic and linear terms, i.e., $V^\mathrm{f}(x) = x^\top P^\mathrm{f}x + \nu^\top x$, it is easy to see that for symmetric $\Lambda$ and $P^\mathrm{f}$ Assumption \ref{as:semidef} is equivalent to the conditions 
\[ P^\mathrm{f} \succeq - \Lambda \quad \mbox{ and } \quad \nu = -v - (2\Lambda + 2P^\mathrm{f})\bar x^\star. \]
Assumption \ref{as:bounded}, in turn, is implied by the condition
\[ P^\mathrm{f} \succ - \Lambda \]
in case $v\ne 0$ and is equivalent to $P^\mathrm{f} \succeq \Lambda$ in the special case that $v=0$. Further conditions for linear quadratic case can be found in the companion paper \cite{ZanG24}.

Given that the optimal control problem from Example \ref{ex:simple_lqr1} is strictly dissipative with storage function $\lambda(x)=-cx^2$ for any $c\in(0,1]$, we conclude that any terminal cost of the form $V^{\mathrm{f}}(x)= ax^2$ for any $a>0$ stabilizes the RH closed loop. As in this example we obtain an unstable closed loop without terminal cost, we see that the condition on $a$ is tight here. The next section illustrates our results with a somewhat more complicated example.

\section{Illustrative example}\label{sec:example}

\begin{figure}[htb]
	\begin{center}
	\includegraphics[width=0.8\linewidth]{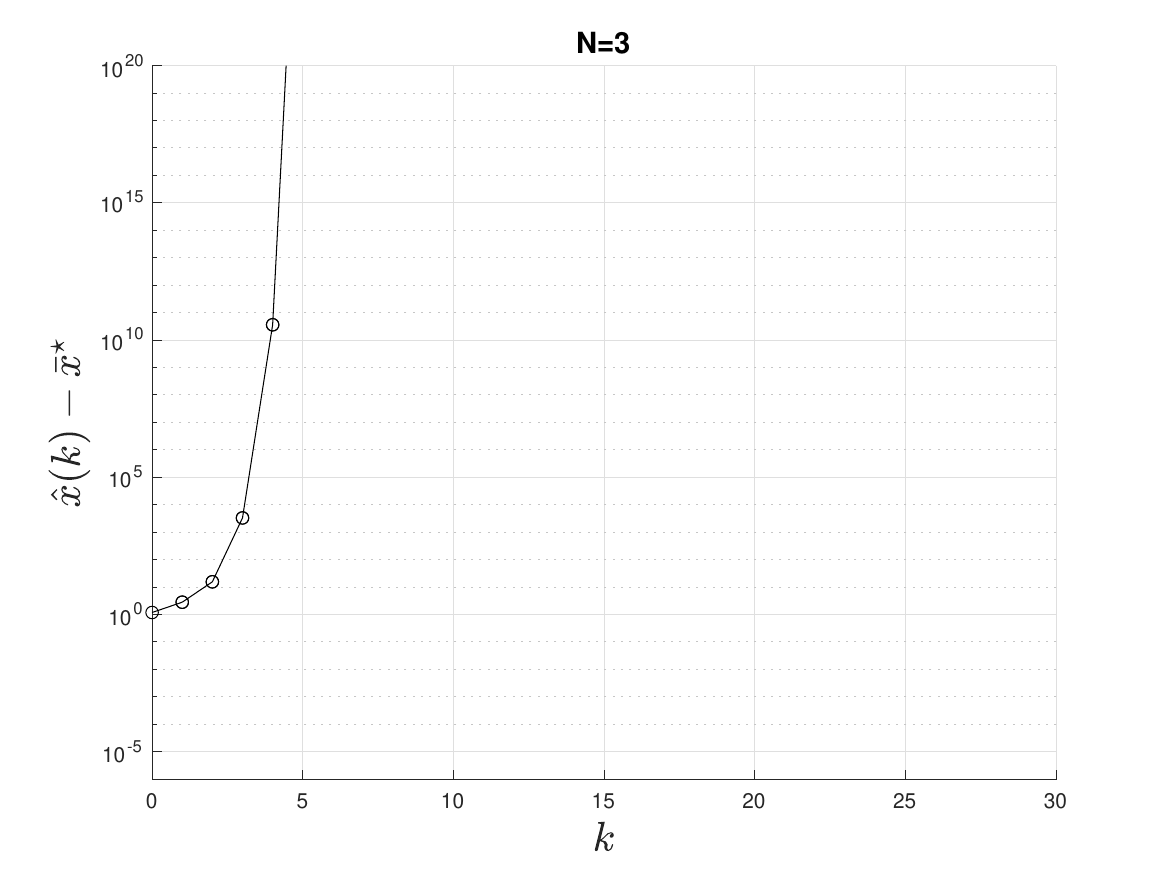}
	\end{center}
	\caption{Difference between closed-loop solution and $\bar x^\star$ with $V^{\mathrm f} \equiv 0$}
	\label{fig:cost0N3}
\end{figure}

We illustrate our findings by an example that can be seen as a more involved nonlinear version of Example \ref{ex:simple_lqr1}. It is given by 
\[ x_{k+1} = 2x_k + x_k^3 + 1 + u_k, \quad \ell(x,u) = u^2\]
with scalar $x$ and $u$. With a little bit of computation one checks that the optimal control problem is strictly dissipative at the equilibrium
\[  \bar x^\star = -\frac{\alpha^2 - 12}{6\alpha}\approx -0.6823278, \quad  \bar u^\star = 0,\]
where $\alpha =(108 + 12\sqrt{93})^{1/3}$, with storage function
\begin{equation} \lambda(x) =  -x^2 + \nu x,\label{eq:lambdaex}\end{equation}
where 
\begin{align*} \nu & = \frac{\alpha(-8\sqrt{93} - 72) + \alpha^2(-6\sqrt{93} - 54) + 144\sqrt{93} + 1392}{\alpha(18\sqrt{93} + 174) - (\alpha^2 - 12)(9 + \sqrt{93})}\\
& \approx -1.3646556.\end{align*}
Figure \ref{fig:cost0N3} shows that the receding horizon closed loop is unstable for terminal cost $V^{\mathrm f} \equiv 0$. The computation was done for $N=3$ but the results are similar for other optimization horizons. 

\begin{figure}[htb]
	\begin{center}
	\includegraphics[width=0.8\linewidth]{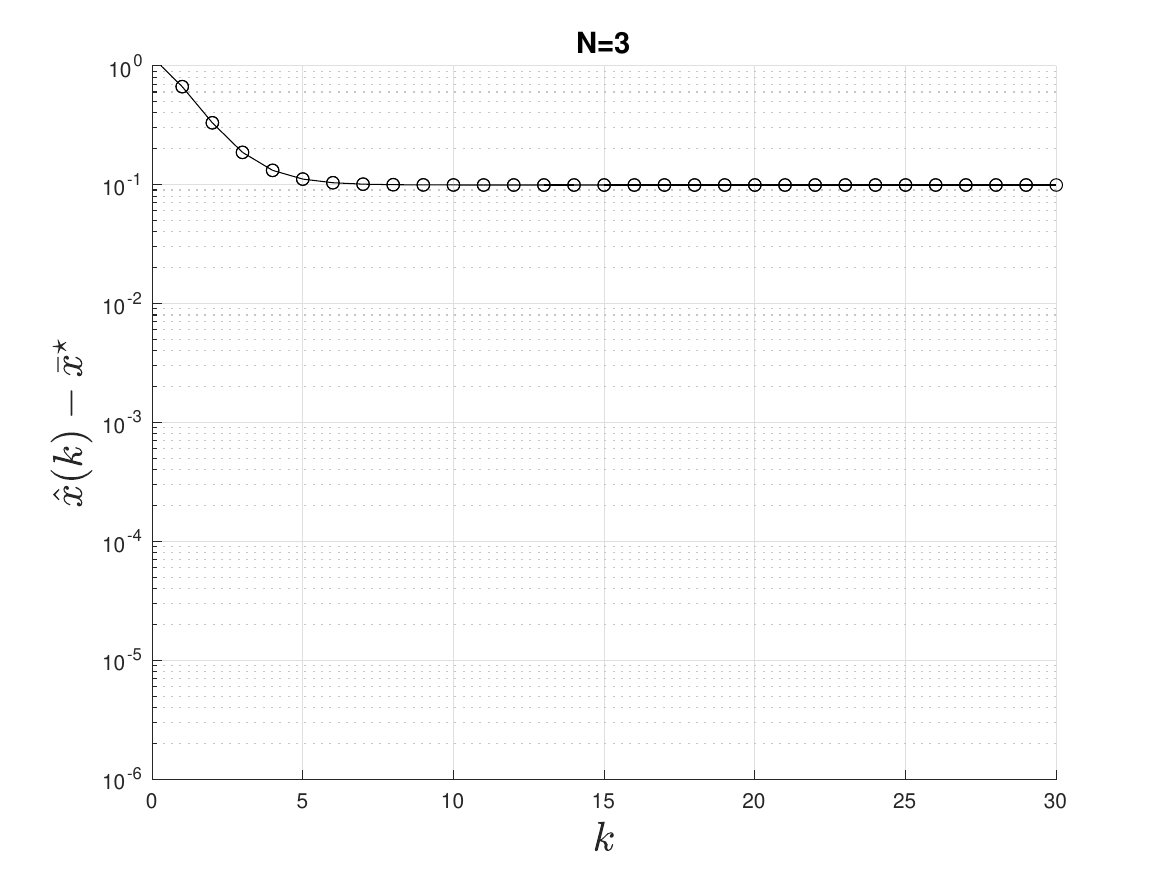}
	\includegraphics[width=0.8\linewidth]{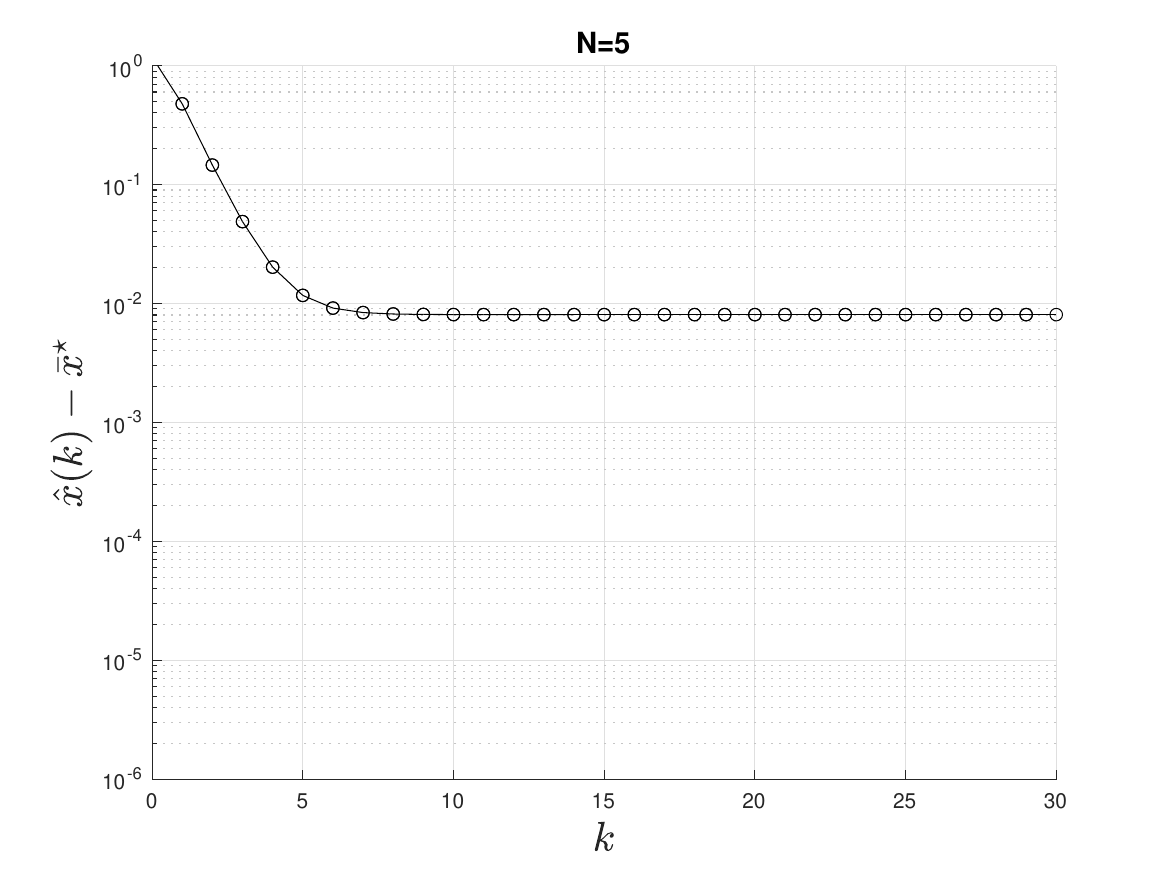}
	\end{center}
	\caption{Difference between closed-loop solution and $\bar x^\star$ with $V^{\mathrm f}(x)=2x^2$ and optimization horizon $N=3$ (top) and $N=5$ (bottom)}
	\label{fig:cost1}
\end{figure}

Figure \ref{fig:cost1} shows the practical asymptotic stability for terminal cost $V^{\mathrm f} = 2x^2$, for which $V^{\mathrm f} + \lambda$ is bounded from below but not positive definite in $\bar x^\star$. The closed-loop solution ends up near the optimal equilibrium $\bar x^\star$ (depicted with the dashed red line), and for the larger value $N=5$ (bottom) it ends up closer to $\bar x^\star$ then for the smaller value $N=3$ (top).

\begin{figure}[htb]
	\begin{center}
	\includegraphics[width=0.8\linewidth]{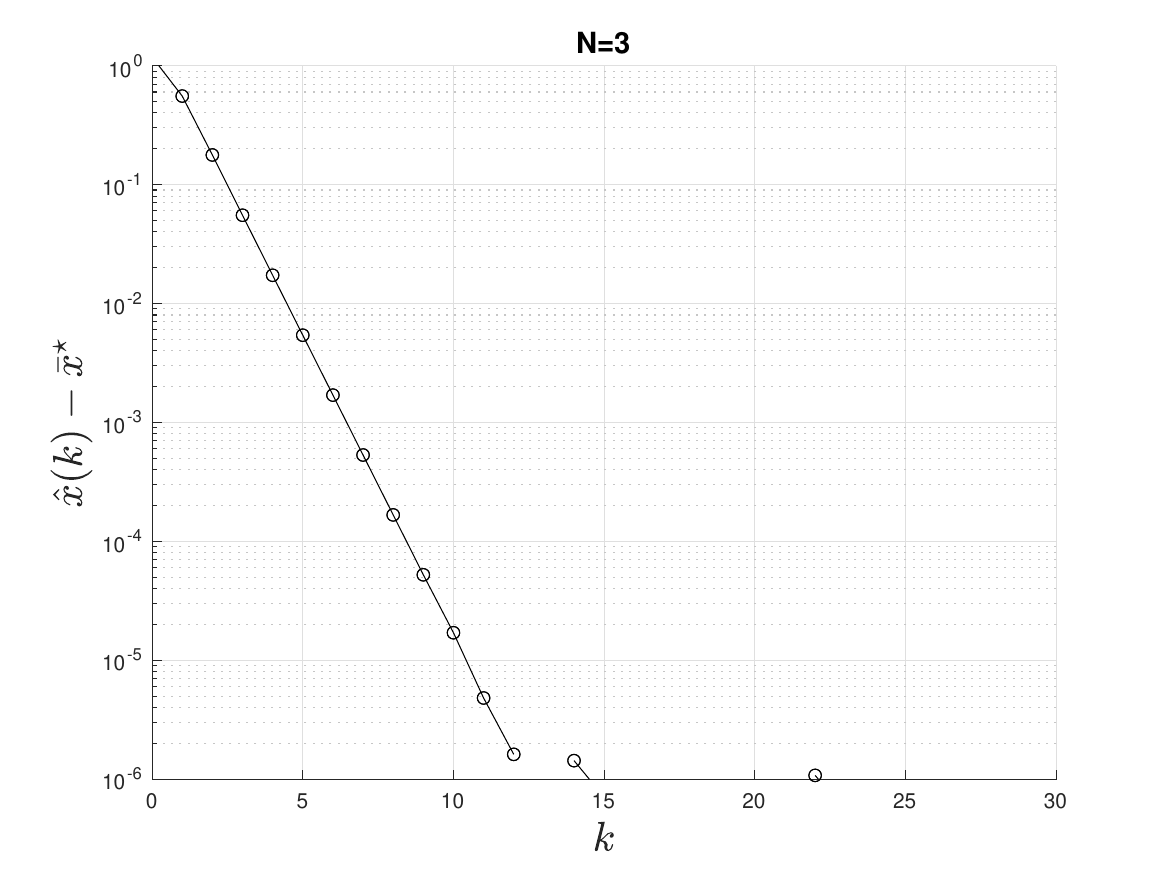}
	\end{center}
	\caption{Difference between closed-loop solution and $\bar x^\star$ with $V^{\mathrm f}(x)=2x^2-2\nu x$ and optimization horizon $N=3$. After time step $k=12$ the solution deteriorates because of roundoff errors.}
	\label{fig:cost2}
\end{figure}

Finally, Figure \ref{fig:cost2} shows that already for optimization horizon $N=3$ ``true'' asymptotic stability (up to roundoff errors) holds for terminal cost $V^{\mathrm f} = x^2-\nu x$, for which $V^{\mathrm f} + \lambda$ vanishes and is thus positive semidefinite at $\bar x^\star$. Note that the inequality on the optimal value function in Theorem \ref{thm:semidef} is satisfied with linear $\rho$, because $V_{\lambda,N}$ grows quadratically and $\inf_u L$ grows at least quadratically\footnote{More precisely, $\inf_u L$ grows quadratically near $\bar x^\star$ and with the order $(x-\bar x^\star)^4$ away from $\bar x^\star$.} in $x-\bar x^\star$, hence the reasoning from Remark \ref{rem:semidef}(ii) applies.

\section{Conclusions}\label{sec:conclusion}

We have shown that receding horizon control with strictly pre-dissipative optimal control problem with storage function $\lambda$ can be stabilized by suitable designed terminal costs $V^{\mathrm{f}}$. For obtaining practical asymptotic stability it is sufficient that $V^{\mathrm{f}}+\lambda$ is bounded from below. If this sum is, in addition, positive definite, then ``true'', i.e., non-practical asymptotic stability can be concluded.

\bibliographystyle{IEEEtranS}
\bibliography{bibliography,bibliography_lars}

\end{document}